\documentclass[11pt]{article}
\usepackage[russian, english]{babel}
\usepackage{amssymb,amsmath,amsthm}
\usepackage{graphicx}
\usepackage{epstopdf}
\usepackage{setspace}
\usepackage{hyperref}
\newcommand{\bc}{\begin{center}}
	\newcommand{\ec}{\end{center}}
\newcommand{\be}{\begin{equation}}
	\newcommand{\ee}{\end{equation}}
\newcommand{\ba}{\begin{array}}
	\newcommand{\ea}{\end{array}}
\newcommand{\edc}{\end{document}}

\begin{document}
\sloppy

\begin{center}
{\bf \large {Gibbs measures for a Hard-Core model with a countable set of states}}
\end{center}
\begin{center}
U.A.Rozikov\footnote{V.I.Romanovsky Institute of Mathematics, Tashkent, Uzbekistan;}$^,$
\footnote{Central Asian University, 264, Milliy bog str. Tashkent, 111221, , Uzbekistan;}$^,$
\footnote{National University of Uzbekistan named after Mirzo Ulugbek, Tashkent, Uzbekistan.},
R.M.Khakimov $^{1,}$\footnote{Namangan State University, Namangan, Uzbekistan.\\
 E-mail: rozikovu@yandex.ru, rustam-7102@rambler.ru, mmtmuxtor93@mail.ru},
M.T.Makhammadaliev $^{1}$
\end{center}
\begin{abstract}
In this paper, we focus on studying non-probability Gibbs measures for a Hard Core (HC) model on a Cayley tree of order $k\geq 2$, where the set of integers $\mathbb Z$ is the set of spin values.

It is well-known that each Gibbs measure, whether it be a gradient or non-probability measure, of this model corresponds to a boundary law. A boundary law can be thought of as an infinite-dimensional vector function defined at the vertices of the Cayley tree, which satisfies a nonlinear functional equation. Furthermore, every normalisable boundary law corresponds to a Gibbs measure. However, a non-normalisable boundary law can define gradient or non-probability Gibbs measures.

In this paper, we investigate the conditions for uniqueness and non-uniqueness of translation-invariant and periodic non-probability Gibbs measures for the HC-model on a Cayley tree of any order $k\geq 2$.

{\bf Mathematics Subject Classifications (2010).} 82B26 (primary);
60K35 (secondary)

\textbf{Key words:} HC-model, configuration, Cayley tree, Gibbs measure, non-probability Gibbs measure, boundary law.

\medskip
\end{abstract}
\begin{center}
\textbf{1. INTRODUCTION}
\end{center}

The Cayley tree $\Im^k=(V, L)$ of order $ k\geq 1 $ is an infinite tree, i.e. graph without cycles, each vertex of which has exactly $k+1$ edges, where $V$ is the set of vertices of $\Im^k$ and $L$ is the set of edges. If  $l \in L$ an edge with  endpoints  $x, y\in V$ then we write $l=\langle x,y\rangle $ and the endpoints are called nearest neighbors.

%

There are a lot of papers devoted to the study of the limiting Gibbs measures for HC models with a finite number of spin values on the Cayley tree (see, for example \cite{RM}, \cite{R}, \cite{BR}), as well as the papers, devoted to the study of (gradient) Gibbs measures for models with an infinite set of spin values (see, for example, \cite{GR}-\cite{HK1}).

In this paper, we study the non-probability Gibbs measures for a Hard Core model with a countable set of spin values on a Cayley tree of order $k\geq 2$. In the paper \cite{RM1} for this model the non-probability Gibbs measures for $k=2, 3, 4$ are studied.
The conditions for the uniqueness and non-uniqueness of translation-invariant (TI) and periodic non-probability Gibbs measures (PNPGM) are found. We note that PNPGMs are already used in physics and other fields of science (see, for example, \cite{AI} and references therein).

In this paper, we generalize the results from \cite{RM1}. The notion of $q$-periodic Gibbs measure is introduced. It is shown that for the model under consideration there is no $q$-periodic Gibbs measure if $q$ is odd. For $q=2$, the uniqueness of the TINPGM and the non-uniqueness of the PNPGMs for the HC-model with a countable set of spin values on the Cayley tree of arbitrary order are proved. In addition, for $q=4$ on a Cayley tree of arbitrary order, the conditions and the exact value of the parameter $\lambda_{cr}^{(1)}$ (resp. $\lambda_{cr}^{(3)}$) are found such that for $0<\lambda\leq\lambda_{cr}^{(1)}$ (resp. $\lambda \geq\lambda_{cr}^{(3)}$) there is exactly one TINPGM (resp. PNPGM), and for $ \lambda>\lambda_{cr}^{(1)}$ (resp. $0<\lambda<\lambda_{cr}^{(3)}$) there are exactly three TINPGMs (resp. PNPGMs).

Let us define our HC model  with a countable set of states on the Cayley tree. The configuration $\sigma = \{\sigma(x) |x \in V \}$ on the Cayley tree is given as a function from $V$ to the set $\mathbb{Z}$, i.e. in this model, each vertex $x$ is assigned one of the values $\sigma(x)\in \mathbb{Z}$.

We consider the set $\mathbb{Z}$ as the set of vertices of an infinite graph $G$. We use the graph $G$, to define a $G$-admissible configuration as follows. A configuration $\sigma$ is called a $G$-\textit{admissible configuration} on a Cayley tree, if $\{\sigma (x),\sigma (y)\}$ is one edge of the graph $G$ for any pair of nearest neighbors $x,y$ in $V$. Denote the set of $G$-admissible configurations by $\Omega^G$.

The activity set \cite{bw} for the graph $G$ is the bounded function $\lambda : G \mapsto \mathbb{R}_+$ ($\mathbb{R}_+ $ is the set of positive real numbers). The value $\lambda_i$ of the function $\lambda$ at the vertex $i \in \mathbb{Z}$ is called its ``activity''.

For fixed $x^0\in V$
write $x\prec y$ if the path from $x^0$ to $y$ goes through $x$. A vertex $y$ is called a direct successor of a vertex $x$ if $y\succ x$ and $x,y$ are nearest neighbors. Note that in $\Im^k$ any vertex $x\neq x^0$ has $k$ direct successors, and $x^0$ has $k+1$ direct successors. The set of direct successors of the vertex $x$ will be denoted by $S(x)$.

For given $G$ and $\lambda$ we define the Hamiltonian of the $G-$HC-model as
\begin{equation}\label{e1} H^{\lambda}_{G}(\sigma)=\left\{%
\begin{array}{ll}
    J \sum\limits_{x\in{V}}{\ln\lambda_{\sigma(x)},} \ \ $ if $ \sigma \in\Omega^G $,$ \\
   +\infty ,\  \  \  \ \ \ \ \ \ \ \ \ \ $  \ if $ \sigma \ \notin \Omega^G $,$ \\
\end{array}
\right.
\end{equation}
where $J\in \mathbb R$.

\begin{figure}
\begin{center}
   \includegraphics[width=13cm]{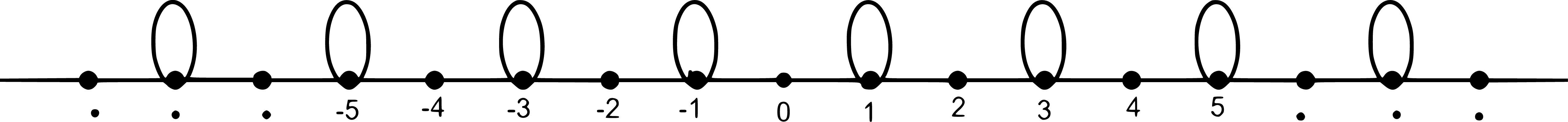}
\end{center}
     \caption{Countable graph $G=\textit{wand}$ with the set of vertices $\mathbb Z$.}\label{fi}
\end{figure}

We consider a specific graph $G=\textit{wand}$ given in Fig. \ref{fi}.
In the case of the graph from Fig. \ref{fi} (see pp. 430-433 in \cite{RM1}), the problem of describing (non-probability) Gibbs measures is reduced to finding solutions to the system of equations
\begin{equation}\label{e2}
\begin{cases} z_{2i+1,x}=\lambda_{2i+1}\prod_{y \in S(x)} \frac{z_{2i,y}+z_{2i+1,y}+z_{2i+2,y}}{z_{-1,y}+z_{1,y}}, \ \ i\in \mathbb{Z}, \\[2 mm]
z_{2i,x}=\lambda_{2i}\prod_{y \in S(x)} \frac{z_{2i-1,y}+z_{2i+1,y}}{z_{-1,y}+z_{1,y}}, \ \ i\in \mathbb{Z}.
\end{cases}
\end{equation}
In \cite{RM1} some solutions of this equation are obtained for $k=2, 3, 4$. Here we consider the case of an arbitrary $k\geq 2$ and generalize the results of \cite{RM1}.

\begin{center}
\textbf{2. TRANSLATION-INVARIANT GIBBS MEASURES}
\end{center}

We consider translation-invariant solutions of the system of equations (\ref{e2}), i.e. $z_x=z\in \mathbb R_{+}^{\infty}.$
In this case, the system of equations (\ref{e2}) has the following form:
\begin{equation}\label{e3}
\begin{cases} z_{2i+1}=\lambda_{2i+1}\Big(\frac{z_{2i}+z_{2i+1}+z_{2i+2}}{z_{-1}+z_{1}}\Big)^k, \ \ i\in \mathbb{Z}, \\[2 mm]
z_{2i}=\lambda_{2i}\Big(\frac{z_{2i-1}+z_{2i+1}}{z_{-1}+z_{1}}\Big)^k, \ \ i\in \mathbb{Z}.
\end{cases}
\end{equation}
Here $\lambda_i>0, \ z_i>0$.

\textbf{Remark 1.} Note that $\lambda_0=1$ and $z_0=1$ (normalization at $0$).

\textbf{Definition 1}. If $z_{i+q}=z_i$ for some $q\geq 1$ and for any $i\in \mathbb{Z}$, then the sequence $\{z_i\}_{i\in \mathbb Z} $ is called $q$-periodic.

The following lemma is true.

\textbf{Lemma 1.} For any $q$, the following statements hold for the system of equations (\ref{e3}):

1. If $z_{i+q}=z_i$, then $\lambda_{i+q}=\lambda_i$.

2. If $q$ is odd, then the system of equations (\ref{e3}) has no $q$-periodic solutions.

\textbf{Proof.}
1. Let $z_{i+q}=z_i$. Then from (\ref{e3}), we have
\begin{equation}\label{e4}
\begin{cases} \lambda_{2i+1}\Big(\frac{z_{2i}+z_{2i+1}+z_{2i+2}}{z_{-1}+z_{1}}\Big)^k=\lambda_{2i+q+1}\Big(\frac{z_{2i+q}+z_{2i+q+1}+z_{2i+q+2}}{z_{-1}+z_{1}}\Big)^k, \ \ i\in \mathbb{Z}, \\[2 mm]
\lambda_{2i}\Big(\frac{z_{2i-1}+z_{2i+1}}{z_{-1}+z_{1}}\Big)^k=\lambda_{2i+q}\Big(\frac{z_{2i+q-1}+z_{2i+q+1}}{z_{-1}+z_{1}}\Big)^k, \ \ i\in \mathbb{Z}.
\end{cases}
\end{equation}
Hence, it is clear that $z_{2i-1}=z_{2i+q-1}$, $z_{2i}=z_{2i+q}$, $z_{2i+1}=z_{2i+q+1}$ and $z_{2i+2}=z_{2i+q+2}$. Then $\lambda_{i+q}=\lambda_i$.

2. Let $z_{i+q}=z_i$. Suppose that $q=2l+1, l\in\mathbb Z$. Then from the system of equations (\ref{e4}) we get:
\begin{equation}\label{e5}
\begin{cases} \lambda_{2i+1}\Big(\frac{z_{2i}+z_{2i+1}+z_{2i+2}}{z_{-1}+z_{1}}\Big)^k=\lambda_{2j}\Big(\frac{z_{2j-1}+z_{2j+1}}{z_{-1}+z_{1}}\Big)^k,\\[2 mm]
\lambda_{2i}\Big(\frac{z_{2i-1}+z_{2i+1}}{z_{-1}+z_{1}}\Big)^k=\lambda_{2j-1}\Big(\frac{z_{2j-2}+z_{2j-1}+z_{2j}}{z_{-1}+z_{1}}\Big)^k,
\end{cases}  i\in \mathbb{Z}, \ j=i+\frac{q+1}2.
\end{equation}
By the first assertion of Lemma 1, we have $\lambda_{i+q}=\lambda_i$. From the first equation of
(\ref{e5}) we get
$$z_{2i}+z_{2i+1}+z_{2i+2}=z_{2i+q}+z_{2i+q+2}.$$
Since $z_{i}=z_{i+q}$, it follows from the last equation that $z_{2i+1}=0$, i.e. the system of equations (\ref{e5}) has no solutions $z_i>0$. Lemma is proved.

\textbf{Corollary 1}. If $z_{i+q}=z_i$ in the system of equations (\ref{e3}), then  $\lambda_i\in\big\{1,\lambda_1,\lambda_2,\cdots,\lambda_{q-1}\big\}$.
In particular,  $\lambda_{2i}=1$ and $\lambda_{2i+1}=\lambda$ for $q=2$.

\textbf{Corollary 2}. The 2-periodic solutions of the system of equations (\ref{e3}) have the form:
$$(\cdots,a^*,1,a^*,1,a^*,1,a^*,1,a^*,1,\cdots).$$

\textbf{Proof.} It is clear that $z_{i+2}=z_i$ holds for the 2-periodic solutions of the system of equations (\ref{e3}). Then taking into account $z_0=1$, we get $z_{2i}=1$. Accordingly,  $z_1=z_{2i+1}$ and $\lambda_1=\lambda_{2i+1}$ for $q=2$.

\textbf{Corollary 3.} The 4-periodic solutions of the system of equations (\ref{e3}) have the form:
$$(\cdots,1,a^*,\lambda_2,c^*,1,a^*,\lambda_2,c^*,1,a^*,\lambda_2,c^*,1,a^*,\lambda_2,c^*,1,\cdots).$$

\textbf{Proof.} The proof is similar to the proof of Corollary 2.

\textbf{Definition 2}. The Gibbs measures corresponding to the $q$-periodic solutions of the system of equations (\ref{e3}) are called $q$-periodic Gibbs measures.

\textbf{Case $q=2$.} By Corollary 2, we have that the 2-periodic solutions (\ref{e3}) have the form
$$
z_j= \begin{cases}
1, \ \ \mbox{if} \ \ j\equiv0\mod2,\\
a, \ \ \mbox{if} \ \ j\equiv1\mod2.
\end{cases}
$$
In this case, due to Corollary 1, we write the system of equations (\ref{e3}) as follows:
\begin{equation}\label{e6}
a=\lambda\Big(\frac{a+2}{2a}\Big)^k.
\end{equation}
Then from (\ref{e6}) we obtain the equation
$$g(a)=2^ka^{k+1}-\lambda(a+2)^k=0,$$
which, by the Descartes rule of signs has at most one positive solution.
 On the other hand, $g(0)=-2^k\lambda<0$ and $g(a)\rightarrow +\infty$ for $a\rightarrow +\infty$, i.e. the equation $g(a)=0$ has at least one positive solution. Therefore, this equation has exactly one positive solution $a^*=a^*(k;\lambda)$ for any $\lambda>0$.

We find the coordinates of the vector $z$, which is the solution of the system of equations (\ref{e3}), as follows:
\begin{equation}\label{e7}
\left(\dots,a^*,1,a^*,1,a^*,1,a^*,1,a^*,1,a^*,\dots\right).
\end{equation}

\textbf{Remark 2.} Obviously, the series obtained from the sequence of solutions (\ref{e7}) diverges. Hence, the measure corresponding to this solution is not normalisable. In addition, the Gibbs measures corresponding to the Hamiltonian of the HC model in our case are not gradient.

Thus, the following theorem is true.

\textbf{Theorem 1.} \emph{Let $k\geq2$ and $q=2$. Then for the HC-model with a countable set of spin values (corresponding to the graph from Fig.1) for any $\lambda>0$, there is exactly one TINPGM $\mu_{0}$ corresponding to a solution of the form (\ref{e7})}.

\textbf{Case $q=4$.} Due to Corollary 3, the 4-periodic solutions of (\ref{e3}) have the form
$$
z_j= \begin{cases}
1, \ \ \mbox{if} \ \ j\equiv0   \mod4, \\
a, \ \ \mbox{if} \ \ j\equiv1 \mod4, \\
\lambda_2, \ \mbox{if} \ \ j\equiv2  \mod4, \\
c, \ \ \mbox{if} \ \ j\equiv3 \mod4.
\end{cases}
$$

In this case, we write the system of equations (\ref{e3}) as follows:
\begin{equation}\label{e8}
\begin{cases}
a=\lambda_1\Big({1+\lambda_2+a \over a+c}\Big)^k, \\
c=\lambda_3\Big({1+\lambda_2+c \over a+c}\Big)^k.
\end{cases}
\end{equation}

The problem of the finding of the general form of solutions of the equation (\ref{e8}) seems
to be very difficult. Therefore, in (\ref{e8}) we assume that $\lambda_{2i+1}=\lambda,$ $i=0,1$ and write the system of equations (\ref{e8}) as follows:
\begin{equation}\label{e9}
\begin{cases}
a=\lambda\Big({1+\lambda_{2}+a \over a+c}\Big)^k, \\
c=\lambda\Big({1+\lambda_{2}+c \over a+c}\Big)^k.
\end{cases}
\end{equation}

In the system of equations (\ref{e9}), we first consider the case $a=c$. In this case, we get the equation
\begin{equation}\label{e10}
a=\lambda\Big(\frac{1+\lambda_{2}+a}{2a}\Big)^k.
\end{equation}

The following lemma is true.

\textbf{Lemma 2.} Let $k\geq2$. Then for any $\lambda>0$ and $\lambda_{2}>0$ the equation (\ref{e10}) has a unique positive solution.

\textbf{Proof.} The equation (\ref{e10}) could be rewritten as
$x=f(x),$ where
$$f(x)=\lambda\Big(\frac{1+\lambda_{2}+x}{2x}\Big)^k.$$
Note that the derivative of the function $f(x)$ is negative
$$f'(x)=-\frac{k\lambda(1+\lambda_{2})}{2x^2}\cdot\Big(\frac{1+\lambda_{2}+x}{2x}\Big)^{k-1}<0,$$
i.e. the function $f(x)$ decreases for all $x>0$. On the other hand, $f(x)\rightarrow +\infty$ for $x\rightarrow 0+0$ and $f(+\infty)= 2^{-k} \lambda$. Hence, the equation (\ref{e10}) has a unique positive solution $a^*=a^*(k;\lambda,\lambda_{2})$ for any $\lambda>0$ and $\lambda_{2} >0$. Lemma is proved.

The following assertion holds.

\textbf{Proposition 1.} Let $k\geq2$, $\lambda_{2}>0$ and $\lambda_{cr}^{(1)}(k,\lambda_{2})=\frac{2^{k}(\lambda_{2}+1)}{(k-1)k^k}$. Then 
\begin{itemize}
	\item 
if  $0<\lambda\leq\lambda_{cr}^{(1)}$, the system of equations (\ref{e9}) has exactly one solution of the form $(a^*, a^*)$,
\item if  $\lambda>\lambda_{cr}^{(1)}$  the system of equations (\ref{e9}) has exactly three solutions of the form $(a^*, a^*)$, $(a,c)$ and $(c,a)$.
\end{itemize}
\textbf{Proof.}
From (\ref{e9}) we get
\begin{equation}\label{s1}
\frac{a}{c}=\left(\frac{\lambda_2+1+a}{\lambda_2+1+c}\right)^k.
\end{equation}
We introduce the notation $\frac{\lambda_2+1+a}{\lambda_2+1+c}=t$, $t>0$.
Then by virtue of (\ref{s1}) ($a=c\cdot t^k$), after some algebra, we get
$$\left(t-1\right)\left(c\cdot(t^{k-1}+t^{k-2}+\dots+t)-\lambda_2 -1\right)=0.$$
Hence $t=1$ or
$$c\cdot(t^{k-1}+t^{k-2}+\dots+t)-\lambda_2 -1=0.$$
It is clear that for $t=1$ we have the solution $a=c=a^*$. By Lemma 2 it follows that for any $\lambda>0$ and $\lambda_{2}>0$ a solution of this kind is unique.

Let $t\neq 1$. Then
$$c(t)=\frac{\lambda_2+1}{t^{k-1}+t^{k-2}+\dots+t}$$
and the function $c(t)$ is uniquely determined for each value of $t$, since
$$c'(t)=-(\lambda_2+1)\frac{(k-1)t^{k-2}+(k-2)t^{k-2}+\dots+1}{(t^{k-1}+t^{k-2}+\dots+t)^2}<0.$$
The values $a(t)$ corresponding to each value $c(t)$ are determined by the formula
$$a(t)= t^k\cdot c(t)=\frac{(\lambda_2+1)t^k}{t^{k-1}+t^{k-2}+\dots+t}.$$
We substitute the expressions for $a(t)=c(t)\cdot t^k$ and $c(t)$ into the first equation in $(\ref{e9})$. Then
\begin{equation}\label{s2}
\frac{(\lambda_2+1)t^k}{t^{k-1}+t^{k-2}+\dots+t}=\lambda\left(t^k+t^{k-1}+t^{k-2}+\dots+t \over t^k+1\right)^k.
\end{equation}
The equation (\ref{s2}) has a solution $t=t(\lambda,\lambda_2)$, but it is very hard to solve. Therefore, we consider the equation (\ref{s2}) with respect to the variable $\lambda$, i.e. from (\ref{s2}) we find that
\begin{equation}\label{s3}
\lambda(t)=\frac{(\lambda_2+1)(t^k+1)^k}{(t^{k-1}+t^{k-2}+\dots+t)(t^{k-1}+t^{k-2}+\dots+1)^k}=
\frac{(\lambda_2+1)(t^k+1)^k}{\Big(\sum_{i=1}^{k-1}t^i\Big)\Big(\sum_{i=0}^{k-1}t^i\Big)^k}.
\end{equation}
Let us prove that each value of $\lambda$ corresponds to only one value of $t$. Note that if $t$ is a solution to (\ref{s3}), then $\frac{1}{t}$ is also a solution to (\ref{s3}). Hence, it suffices to show that each value of $\lambda$ corresponds to exactly one value $t>1$ (or $t<1$).
To do this, consider the derivative of the function $\lambda(t)$:
$$\lambda'(t)=\frac{(\lambda_2+1)(t^k+1)^{k-1}\cdot\vartheta(t,k)}{t^2(t^k-1)(t^{k-1}-1)^2\cdot(t^{k-1}+t^{k-2}+\dots+t+1)^k},$$
where
$$\vartheta(t,k)=(t^k-1)(t^{2k}-1)+kt(t^{k-2}-1)(t^{2k}-1)-2k^2t^k(t-1)(t^{k-1}-1).$$
It can be seen that $t=1$ is a double root of the polynomial $\vartheta(t,k)$. Let us show that $t=1$ is a root of $\vartheta(t,k)$ multiplicity four. To do this, we introduce the notation $t-1=x$ ($t=x+1$) and prove that $x=0$ is the root of $\vartheta_1(x,k)=\vartheta_1(t-1,k)$ multiplicity four, i.e. in the polynomial $\vartheta_1(x,k)$ the least power of the variable $x$ is four.
 To do this, we introduce notation 
 $$C_k^i={k\choose i}$$ and rewrite $\vartheta_1(x,k)$ as follows:
$$\vartheta_1(x,k)=\left(\sum_{i=1}^{2k}{C_{2k}^ix^i}\right)\left(\sum_{i=1}^{k}{C_{k}^ix^i}+k(x+1)\sum_{i=1}^{k-2}{C_{k-2}^ix^i}\right)
-2k^2\left(\sum_{i=0}^{k}{C_{k}^ix^{i+1}}\right)\left(\sum_{i=1}^{k-1}{C_{k-1}^ix^i}\right).$$
It is easy to show that the coefficients for $x^2$ and $x^3$ in the polynomial $\vartheta_1(x,k)$ are equal to zero, while for $x^4$ it is different from zero. Indeed, the coefficient at $x^2$:
$$C_{2k}^1\cdot\Big(C_k^1+kC_{k-2}^1\Big)-2k^2C_k^0\cdot C_{k-1}^1=2k\big(k+k(k-2)\big)-2k^2(k-1)=0,$$
the coefficient at $x^3$:
$$C_{2k}^1\cdot\Big(C_k^2+kC_{k-2}^2+kC_{k-2}^1\Big)+C_{2k}^2\Big(C_k^1+kC_{k-2}^1\Big)-2k^2\Big(C_k^0\cdot C_{k-1}^2+C_k^1\cdot C_{k-1}^1\Big)=$$
$$=k^2(3k^2-5k+2)-k^2(3k^2-5k+2)=0,$$
the coefficient at $x^4$:
$$C_{2k}^1\cdot\Big(C_k^3+kC_{k-2}^3+kC_{k-2}^2\Big)+
C_{2k}^2\Big(C_k^2+kC_{k-2}^2+kC_{k-2}^1\Big)+C_{2k}^3\Big(C_k^1+kC_{k-2}^1\Big)-$$
$$-2k^2\Big(C_k^0\cdot C_{k-1}^3+C_k^1\cdot C_{k-2}^2+C_k^2\cdot C_{k-2}^1\Big)=\frac{k^2(k-1)(2k^2-k+3)}6\neq0.$$

Next, we will show that $\vartheta_1(x,k)\geq0$, i.e.
\begin{equation}\label{s4}
\left(\sum_{i=1}^{2k}{C_{2k}^ix^i}\right)\left(\sum_{i=1}^{k}{C_{k}^ix^i}+k(x+1)\sum_{i=1}^{k-2}{C_{k-2}^ix^i}\right)
>\left(2k\sum_{i=0}^{k}{C_{k}^ix^{i+1}}\right)\left(k\sum_{i=1}^{k-1}{C_{k-1}^ix^i}\right).
\end{equation}
Let us first show that in (\ref{s4}) the expression in the second factor of the LHS is greater than the expressions in the second factor of the RHS, i.e. validity of the inequality
$$\sum_{i=1}^{k}{C_{k}^ix^i}+k(x+1)\sum_{i=1}^{k-2}{C_{k-2}^ix^i}>k\sum_{i=1}^{k-1}{C_{k-1}^ix^i}.$$
Indeed,
$$\sum_{i=1}^{k}{C_{k}^ix^i}+k(x+1)\sum_{i=1}^{k-2}{C_{k-2}^ix^i}=\sum_{i=1}^{k}{C_{k}^ix^i}+k\sum_{i=1}^{k-2}{C_{k-2}^ix^i}+k\sum_{i=2}^{k-1}{C_{k-2}^{i-1}x^i}=$$
$$=\sum_{i=1}^{k}{C_{k}^ix^i}+ k\Big((k-2)x+\sum_{i=2}^{k-2}\Big(C_{k-2}^i+C_{k-2}^{i-1}\Big)x^i+x^{k-1}\Big)=$$
$$=\sum_{i=2}^{k}{C_{k}^ix^i}+k\Big((k-1)x+\sum_{i=2}^{k-1}C_{k-1}^ix^i\Big)=\sum_{i=2}^{k}{C_{k}^ix^i}+k\sum_{i=1}^{k-1}C_{k-1}^ix^i>k\sum_{i=1}^{k-1}C_{k-1}^ix^i$$

Now we show that in the inequality (\ref{s4}) the expression in the first factor of the LHS is greater than the expressions in the first factor of the RHS, i.e. validity of the following inequality:
$$\sum_{i=1}^{2k}{C_{2k}^ix^i}\geq 2k\sum_{i=0}^{k}{C_{k}^ix^{i+1}} = 2k\sum_{i=1}^{k+1}{C_{k}^{i-1}x^i}.$$
To do this, we prove that $C_{2k}^i\geq 2k C_{k}^{i-1}$. Indeed,
$$\frac{(2k)!}{i!(2k-i)!}\geq\frac{2k \cdot  k!}{(i-1)!(k+1-i)!} $$
or
\begin{equation}\label{s5}
(2k-1)!\cdot (k+1-i)! \geq i\cdot k!\cdot(2k-i)!.
\end{equation}
Using mathematical induction, we show the validity of the inequality (\ref{s5}) for $i\geq3$. For $i=3$, it is trivial:
$$(2k-1)!\cdot (k-2)! \geq 3\cdot k!\cdot(2k-3)! \ \ \Rightarrow \ \ k^2-3k+2\geq0 \ \ \Rightarrow \ \ k\geq2. $$
Assume that the inequality (\ref{s5}) is true for $i=p$,
and let us prove that (\ref{s5}) is also true for $i=p+1$:
$$(2k-1)!\cdot (k-p)! \geq (p+1)\cdot k!\cdot(2k-p-1)!$$
Using the inequality for $i=p$, we can get
$$(p+1)\cdot k!\cdot(2k-p-1)!\leq\frac{(p+1)(2k-1)!\cdot (k+1-p)!}{p(2k-p)}.$$
Hence, it suffices to show that
$$\frac{(p+1)(2k-1)!\cdot (k+1-p)!}{p(2k-p)}\leq(2k-1)!\cdot (k-p)!.$$
From the last inequality, after some algebra, we obtain
$$(p+1) (k+1-p)\leq p(2k-p) \ \ \Rightarrow \ \ k(p-1)\geq0.$$
Thus, for $i\geq3$ and $k\geq2$ we have
$$\sum_{i=1}^{2k}{C_{2k}^ix^i}\geq 2k\sum_{i=0}^{k}{C_{k}^ix^{i+1}},$$
i.e. for $k\geq2$
$$\left(\sum_{i=1}^{2k}{C_{2k}^ix^i}\right)\left(\sum_{i=1}^{k}{C_{k}^ix^i}+k(x+1)\sum_{i=1}^{k-2}{C_{k-2}^ix^i}\right)
>\left(2k\sum_{i=0}^{k}{C_{k}^ix^{i+1}}\right)\left(k\sum_{i=1}^{k-1}{C_{k-1}^ix^i}\right).$$
This means that all the coefficients of $x^n$, $n\geq4$, in the polynomial $\vartheta_1(x,k)$ are positive.

It follows that
$$\vartheta(t,k)=(t-1)^4\cdot\eta(t),$$
and $\eta(t)>0$ for $t>0$. It means that
$$\lambda'(t)=\frac{(\lambda_2+1)(t-1)^4(t^k+1)^{k-1}\cdot\eta(t)}
{t^2(t^k-1)(t^{k-1}-1)^2\cdot(t^{k-1}+t^{k-2}+\dots+t+1)^k}.$$
Therefore, the function $\lambda(t)$ decreases for $t<1$, increases for $t>1$, and reaches its minimum for $t=1$ (see Fig.2):
$$\lambda_{\min}(t)=\lambda(1)=\lambda_{cr}^{(1)}(k,\lambda_{2})=\frac{2^{k}(\lambda_{2}+1)}{(k-1)k^k}.$$

\begin{center}
\includegraphics[width=8cm]{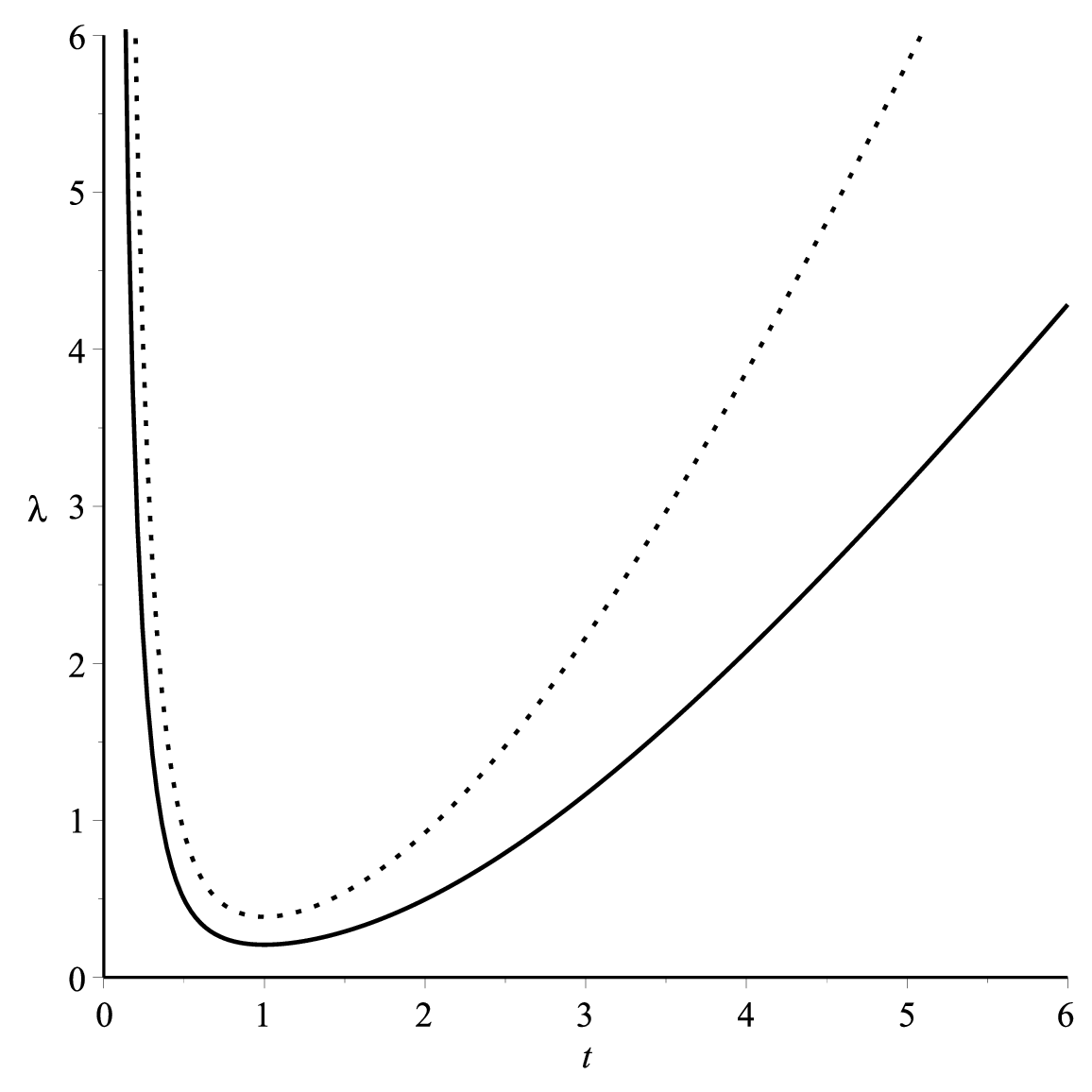}
\end{center}
\begin{center}{\footnotesize \noindent
 Fig.2. Graph of the function $\lambda(t)$ for $k=3$, $\lambda_2=0.4$ (continuous) and $\lambda_2=1.6$ (dotted).}\
\end{center}

Hence, each value $\lambda$ corresponds to only one value $t>1$ (or $t<1$) for $\lambda>\lambda_{cr}^{(1)}$, the value $t=1$ for $\lambda=\lambda_{cr}^{(1)}$ and the equation (\ref{s2}) has no solutions for $\lambda<\lambda_{cr}^{(1)}$. Due to symmetry, the system of equations (\ref{e9}) has a unique solution of the form $(a^*,a^*)$ for $\lambda\leq\lambda_{cr}^{(1)}$, and for $\lambda>\lambda_{cr}^{(1)}$ has exactly three solutions $(a^*,a^*)$, $(a,c)$ and $(c,a)$, $a \neq c$, where $a^*$ is the only positive solution of (\ref{e10}). The proof is complete.

Next, we find the coordinates of the vector $z$, which is the solution of the system of equations (\ref{e3}), as follows:

\begin{itemize}
\item for $0<\lambda\leq\lambda_{cr}^{(1)}$ we get
\begin{equation}\label{e16}
\left(\dots,1,a^*,1,a^*,1,a^*,1,a^*,1,a^*,\dots\right),
\end{equation}
\item for $\lambda>\lambda_{cr}^{(1)}$ we have three solutions
\begin{equation}\label{e17}
\begin{cases}
\left(\dots,1,a^*,\lambda_{2},a^*,1,a^*,\lambda_{2},a^*,1,\dots\right),\\
\left(\dots,1,a,\lambda_{2},c,1,a,\lambda_{2},c,1,a,\lambda_{2},c,\dots\right),\\
\left(\dots,1,c,\lambda_{2},a,1,c,\lambda_{2},a,1,c,\lambda_{2},a,\dots\right).
\end{cases}
\end{equation}
where $a^*$ is the only solution of the equation (\ref{e10}) and $(a,c)$, $(c,a)$ are solutions of the system of equations (\ref{e9}).
\end{itemize}

By Proposition 1, we have

\textbf{Theorem 2.} \emph{Let $k\geq2$, $q=4$ and $\lambda_{cr}^{(1)}(k,\lambda_{2})=\frac{2^ {k}(\lambda_{2}+1)}{(k-1)k^k}$. Then for an HC model with a countable set of spin values (corresponding to the graph in Fig. 1) for $0<\lambda\leq\lambda_{cr}^{(1)}$, there is exactly one TINPGM $\mu_{0}$ corresponding to the solution (\ref{e16}), for $\lambda>\lambda_{cr}^{(1)}$ there are exactly three TINPGMs $\mu_{0}$, $\mu_{1}$, $\mu_{ 2}$ corresponding to solutions (\ref{e17}).}

\textbf{Remark 3.} In \cite{RM1} the result of Theorem 2 was obtained in the cases $\lambda_2=1$ and $k=2,3,4$. In these cases, the critical values of $\lambda$ found in \cite{RM1} coincide with the values of $\lambda_{cr}^{(1)}(k,\lambda_{2})$, but our proof differs from that of \cite{RM1} for $k=2,3,4$. \\

\begin{center}
\textbf{3. PERIODIC GIBBS MEASURES}
\end{center}

In this section, we generalize the results of \cite{RM1} for the 2-periodic Gibbs measures. Such measures correspond to solutions of the following system of equations (see system (23) in \cite{RM1}):
\begin{equation}\label{e18}
\begin{cases}
z_{2i+1}=\lambda_{2i+1}\Big(\frac{\widetilde{z}_{2i}+\widetilde{z}_{2i+1}+\widetilde{z}_{2i+2}}{\widetilde{z}_{-1}+\widetilde{z}_{1}}\Big)^k, \ \ i\in \mathbb{Z}, \\
\widetilde{z}_{2i+1}=\lambda_{2i+1}\Big(\frac{z_{2i}+z_{2i+1}+z_{2i+2}}{z_{-1}+z_{1}}\Big)^k, \ \ i\in \mathbb{Z}, \\
z_{2i}=\lambda_{2i}\Big(\frac{\widetilde{z}_{2i-1}+\widetilde{z}_{2i+1}}{\widetilde{z}_{-1}+\widetilde{z}_{1}}\Big)^k, \ \ i\in \mathbb{Z}, \\
\widetilde{z}_{2i}=\lambda_{2i}\Big(\frac{z_{2i-1}+z_{2i+1}}{z_{-1}+z_{1}}\Big)^k, \ \ i\in \mathbb{Z}.
\end{cases}
\end{equation}

Let the coordinates of the vector $(z, \widetilde{z})$ be the solution of the system of equations (\ref{e18}), \\
where $z=(\ldots,\ z_{-1},z_{0}, \ z_{1},\ldots)$ and  $\widetilde{z}=(\ldots,\ \widetilde{z}_{-1},\widetilde{z}_{0}, \ \widetilde{z}_{1},\ldots).$

\textbf{Definition 3}. The set of ordered pairs $(z, \widetilde{z})$ obtained from the sequences $\{z_i\}_{i\in \mathbb Z}$ and $\{\widetilde {z_i}\}_{i\in \mathbb Z}$ is called $q$-periodic if
$$(z_{i+q},\widetilde{z}_{i+q})=(z_i,\widetilde{z_i}) \ \ \mbox{for} \ \ q\geq 1 \ \ \mbox{and for any} \ \ i\in \mathbb Z.$$

The following lemma is true.

\textbf{Lemma 3.} For any $q$, the following statements hold for the system of equations (\ref{e18}):

1. If $(z_{i+q},\widetilde{z}_{i+q})=(z_i,\widetilde{z_i})$ , then $\lambda_{i+q}=\lambda_i$.

2. If $q$ is odd, then the system of equations (\ref{e18}) has no $q$-periodic solutions.

\textbf{Proof.}
1. Let $(z_{i+q},\widetilde{z}_{i+q})=(z_i,\widetilde{z_i})$. Then from (\ref{e18})
$$
\begin{cases} \lambda_{2i+1}\Big(\frac{\widetilde{z}_{2i}+\widetilde{z}_{2i+1}+\widetilde{z}_{2i+2}}{\widetilde{z}_{-1}+\widetilde{z}_{1}}\Big)^k=
\lambda_{2i+q+1}\Big(\frac{\widetilde{z}_{2i+q}+\widetilde{z}_{2i+q+1}+\widetilde{z}_{2i+q+2}}{\widetilde{z}_{-1}+\widetilde{z}_{1}}\Big)^k, \ \ i\in \mathbb{Z}, \\
\lambda_{2i+1}\Big(\frac{z_{2i}+z_{2i+1}+z_{2i+2}}{z_{-1}+z_{1}}\Big)^k=
\lambda_{2i+q+1}\Big(\frac{z_{2i+q}+z_{2i+q+1}+z_{2i+q+2}}{z_{-1}+z_{1}}\Big)^k, \ \ i\in \mathbb{Z}, \\
\lambda_{2i}\Big(\frac{\widetilde{z}_{2i-1}+\widetilde{z}_{2i+1}}{\widetilde{z}_{-1}+\widetilde{z}_{1}}\Big)^k=
\lambda_{2i+q}\Big(\frac{\widetilde{z}_{2i+q-1}+\widetilde{z}_{2i+q+1}}{\widetilde{z}_{-1}+\widetilde{z}_{1}}\Big)^k, \ \ i\in \mathbb{Z}, \\
\lambda_{2i}\Big(\frac{z_{2i-1}+z_{2i+1}}{z_{-1}+z_{1}}\Big)^k=
\lambda_{2i+q}\Big(\frac{z_{2i+q-1}+z_{2i+q+1}}{z_{-1}+z_{1}}\Big)^k, \ \ i\in \mathbb{Z}.
\end{cases}
$$
Hence, it is clear that $z_{i+q}=z_{i}$ and $\widetilde{z}_{i+q}=\widetilde{z}_{i}$. Then $\lambda_{i+q}=\lambda_i$.

2. Let $(z_{i+q},\widetilde{z}_{i+q})=(z_i,\widetilde{z_i})$. We assume that $q=2j+1, j\in\mathbb Z$. Then from the system of equations (\ref{e18}) we get:
\begin{equation}\label{e19}
\begin{cases} \lambda_{2i+1}\Big(\frac{\widetilde{z}_{2i}+\widetilde{z}_{2i+1}+\widetilde{z}_{2i+2}}{\widetilde{z}_{-1}+\widetilde{z}_{1}}\Big)^k=
\lambda_{2j}\Big(\frac{\widetilde{z}_{2j-1}+\widetilde{z}_{2j+1}}{\widetilde{z}_{-1}+\widetilde{z}_{1}}\Big)^k,\\[2 mm]
\lambda_{2i+1}\Big(\frac{z_{2i}+z_{2i+1}+z_{2i+2}}{z_{-1}+z_{1}}\Big)^k=
\lambda_{2j}\Big(\frac{z_{2j-1}+z_{2j+1}}{z_{-1}+z_{1}}\Big)^k,\\[2 mm]
\lambda_{2i}\Big(\frac{\widetilde{z}_{2i-1}+\widetilde{z}_{2i+1}}{\widetilde{z}_{-1}+\widetilde{z}_{1}}\Big)^k=
\lambda_{2j-1}\Big(\frac{\widetilde{z}_{2j-2}+\widetilde{z}_{2j-1}+\widetilde{z}_{2j}}{\widetilde{z}_{-1}+\widetilde{z}_{1}}\Big)^k,\\[2 mm]
\lambda_{2i}\Big(\frac{z_{2i-1}+z_{2i+1}}{z_{-1}+z_{1}}\Big)^k=\lambda_{2j-1}\Big(\frac{z_{2j-2}+z_{2j-1}+z_{2j}}{z_{-1}+z_{1}}\Big)^k.
\end{cases} i\in \mathbb{Z}, \ j=i+\frac{q+1}2.
\end{equation}
By the first assertion of Lemma 3, we have $\lambda_{2i+1}=\lambda_{2i+q+1}$. From the second equation (\ref{e19}) we get
$$z_{2i}+z_{2i+1}+z_{2i+2}=z_{2i+q}+z_{2i+q+2}.$$
Since $z_{i}=z_{i+q}$, it follows from the last equation that $z_{2i+1}=0$, i.e. the system of equations (\ref{e19}) has no solutions $z_i>0$. Lemma is proved.

\textbf{Corollary 4.} If the system of equations (\ref{e19}) satisfies $(z_{i+q},\widetilde{z}_{i+q})=(z_i,\widetilde{z_i})$, then   $\lambda_i\in\big\{1,\lambda_1,\lambda_2,\cdots,\lambda_{q-1}\big\}$, because $\lambda_0=1$ (normalization at $0$). In particular,  $\lambda_{2i}=1$ and $\lambda_{2i+1}=\lambda$ for $q=2$.

\textbf{Remark 4}. Similarly to the translation-invariant case, for $i=0$ we obtain $\lambda_0=1$ and $z_0=1$ $(\widetilde{z}_0=1)$.

\textbf{Corollary 5.} The 2-periodic solutions of the system of equations (\ref{e18}) have the form $(A,\widetilde{A}),$ where
$$A=(\cdots,a^*,1,a^*,1,a^*,1,a^*,1,\cdots), \ \ \widetilde{A}=(\cdots,\widetilde{a},1,\widetilde{a},1,\widetilde{a},1,\widetilde{a},1,\cdots).$$

\textbf{Proof.}  It is clear that for the 2-periodic solutions of the system of equations (\ref{e18}) satisfy $(z_{i+2},\widetilde{z}_{i+2})=(z_i,\widetilde{z}_i)$. Then taking into account $(z_0,\widetilde{z}_0)=(1,1)$ (normalization at $0$), we get $(z_{2i},\widetilde{z}_{2i})=(1,1 )$. Accordingly, for $q=2$, $(z_{2i+1},\widetilde{z}_{2i+1})=(z_1,\widetilde{z}_1)$ and $\lambda_1=\lambda_{ 2i+1}$.

\textbf{Corollary 6.} The 4-periodic solutions of the system of equations (\ref{e18}) have the form $(A,\widetilde{A}),$
where
$$A=(\cdots,c^*,1,a^*,\lambda_2,c^*,1,a^*,\lambda_2,c^*,1,a^*,\cdots), \ \ \widetilde{A}=(\cdots,\widetilde{c},1,\widetilde{a},\lambda_2,\widetilde{c},1,\widetilde{a},\lambda_2,\widetilde{c},1,\widetilde{a},\cdots).$$

\textbf{Proof}. The proof is similar to the proof of Corollary 5.

\textbf{Case $q=2$}. By virtue of Corollaries 4 and 5, we have that the 2-periodic solutions (\ref{e18}) have the form:
\begin{center}
$z_j= \begin{cases}
1, \ \ \mbox{if} \ \ j\equiv0\mod2, \\
a, \ \ \mbox{if} \ \ j\equiv1\mod2,
\end{cases}$ \ \ \ \
$\widetilde{z}_j= \begin{cases}
1, \ \ \mbox{if} \ \ j\equiv0\mod2, \\
c, \ \ \mbox{if} \ \ j\equiv1 \mod2,
\end{cases}$
\end{center}
In this case, we write the system of equations (\ref{e18}) as follows:
\begin{equation}\label{e20}
\begin{cases}
a=f(c), \\
c=f(a),
\end{cases}
\end{equation}
where
$$f(x)=\lambda\left(\frac{2+x}{2x}\right)^{k}.$$

Before analyzing the system of equations (\ref{e20}), we introduce the following definition.

\textbf{Definition 4.} \cite{Kel} A twice continuously differentiable function $h:[0, \infty) \mapsto [0, \infty)$ is called $S$-shaped if it has the following properties:

$(1)$ function $h(x)$ increasing on $[0, \infty)$ with $h(0)>0$ and $\sup_x{h (x)}<\infty$;

$(2)$ there exists $\widetilde{x}\in (0,\infty)$ such that the derivative $h'$ is monotone increasing in the interval $(0, \widetilde{x})$ and monotone decreasing in the interval $(\widetilde{x },\infty)$; in other words, $\widetilde{x}$ satisfies $h''(\widetilde{x})=0$ and is the only inflection point of $h(x)$.

It is known \cite{Kel} that an $S$-shaped function has at most three fixed points in the interval $[0, \infty)$.

\textbf{Lemma 4}. (Kesten) \cite{K} Let $f:[0,1]\rightarrow [0,1]$ be a continuous function with fixed point $\xi \in (0,1)$. Assume that $f$ is differentiable at the point $\xi$ and $f^{'}(\xi)<-1.$ Then there are points $x_0, \ x_1, \ 0\leq x_0<\xi<x_1 \leq1 $
such that $f(x_0)=x_1$ and $f(x_1)=x_0.$

\textbf{Remark 5}. Note that Kesten's lemma is true if $f:[\alpha, \beta]\rightarrow [\alpha, \beta]$. In our case, the function $f(x)$ is decreasing and $\alpha=\frac{\lambda}{2^k}$, $\beta=\lambda\left(\frac{2^k}{\lambda}+\frac12\right)^k$.

The following proposition is true.

\textbf{Proposition 2}. Let $k\geq2$ and $\lambda_{cr}^{(2)}(k)=\frac{2^{k+1}(k-1)^{k+1}}{k^k} $. Then the system of equations (\ref{e20}) with $\lambda\geq\lambda_{cr}^{(2)}(k)$ has exactly one solution of the form $\left(a^{\ast},a^{ \ast}\right)$, and for $0<\lambda<\lambda_{cr}^{(2)}(k)$ has exactly three solutions of the form $\left(a^{\ast},a^{\ast}\right),$ $\left(a_1,a_2\right),$ $\left(a_2,a_1\right)$.

\textbf{Proof}. Let $h(x)=f(f(x))$. Let us show that the function $h(x)$ is $S$-shaped.
Indeed, the function $h(x)$ for $x>0$ is increasing and bounded, since
$$h'(x)=\frac{k^2\lambda\left(2+f(x)\right)^{k-1}}{2^{k-2}x(x+2)f^k(x)}>0, \ \  h(0)=\frac{\lambda}{2^k}>0,\ \ \lim_{x\rightarrow\infty}{h(x)}=\lambda\left(\frac{2^{k+1}+\lambda}{2\lambda}\right)^k<\infty.$$
Moreover,
$$h''(x)=\frac{\lambda k^2\Big(\frac{f(x)}{2}(k-x-1)+k^2-x-1\Big)\Big(2+f(x)\Big)^{k-2}}{2^{k-4}x^2(x+2)^2f^k(x)}.$$

Let us show that the equation $h''(x)=0$ has exactly one solution. To do this, consider the function
$$\delta(x)=\frac{\lambda}{2}\left(\frac{2+x}{2x}\right)^k(k-x-1)+k^2-x-1.$$

Calculate the derivative $\delta'(x)$:
$$\delta'(x)=-\frac{f(x)}{2x^2+4x}\Big(x^2-2(k-1)x+2k(k-1)\Big)-1,$$
where $f(x)=\lambda\left(\frac{2+x}{2x}\right)^k.$
It is easy to see that $x^2-2(k-1)x+2k(k-1)>0$ for $k\geq2$.

Hence $\delta'(x)<0$, i.e. the $\delta(x)$ is decreasing function.
On the other hand, for $k\geq2$
$$\delta(1)=\frac{\lambda}{2}\Big(\frac{3}{2}\Big)^k(k-2)+k^2-2>0, \ \
\delta(k^2)=-\frac{\lambda}2\Big(\frac{2+k^2}{2k^2}\Big)^k(k^2-k+1)-1<0.$$
Therefore, the equation $\delta(x)=0$ and hence the equation $h''(x)=0$ has a unique root for $x>0$. Hence, the function $h(x)$ is $S$-shaped.

In \cite{Kel}, some properties of an $S$-shaped function are described. In particular, the function $h(x)$ has the following properties:

\begin{itemize}
	\item There is $\lambda_{cr}>0$ such that if $\lambda\geq \lambda_{cr}$ then $h'(x)\leq1$ for any $x\geq0$ and $x_0$ is the only fixed point for $h(x)$.
\item If $\lambda<\lambda_{cr},$ then $h(x)$ has three fixed points $x_1<x_0<x_2$, where $f(x_1)=x_2$ and $f(x_2)=x_1$. Moreover, $h'(x_0)>1,$ $h'(x)<1$ for $x\in [0;x_1]\cup[x_2;\infty)$ and the three fixed points converge to $x_0( \lambda_{cr})$ for $\lambda\rightarrow\lambda_{cr}.$
\end{itemize}

It is easy to see that since $h'(x_0)=(f'(x_0))^2$, then $x_0$ is the only fixed point of the function $h(x)$ if and only if $f'(x_0)\geq-1$.

Let $x_0$ be the unique solution of the equation $f(x)=x$. Calculate the derivative $f'(x_0)$:
$$f'(x)=-\frac{\lambda k(2+x)^{k-1}}{2^{k-1}x^{k+1}} \ \ \Rightarrow \ \ f'(x_0)=-\frac{2k}{2+x_0}.$$
After solving the inequality $f'(x_0)<-1$, we obtain $x_0<2k-2$. Then from $x_0=\lambda\left(\frac{2+x_0}{2x_0}\right)^k$ we deduce that the equation $h(x)=x$ has only one fixed point when
$$\lambda\geq\frac{2^{k+1}(k-1)^{k+1}}{k^k}=\lambda_{cr}^{(2)}(k).$$
From the inequality $f'(x_0)<-1$ we get $\lambda<\lambda_{cr}^{(2)}$. Then by Kesten's lemma it follows that the function $h(x)$ has at least three fixed points. On the other hand, under this condition, by the property of an $S$-shaped function, we have at most three fixed points for $h(x)$. Therefore, for $0<\lambda<\lambda_{cr}^{(2)}$ there are exactly three fixed points of the equation $h(x)=x$.

Thus, the system of equations (\ref{e20}) with $\lambda\geq\lambda_{cr}^{(2)}$ has a unique solution of the form $(a^{\ast},a^{\ast})$, and for $\lambda<\lambda_{cr}^{(2)}$, except for the solution $(a^{\ast},a^{\ast})$, it has two positive solutions $\Big(a_1 ,a_2\Big)$ and $\Big(a_2,a_1\Big)$. The proposition is proved.

Since the series obtained from the sequence of solutions diverge, and by virtue of Proposition 2, we have the following result:

\textbf{Theorem 3.} \textit{Let $k\geq2$, $q=2$ and $\lambda_{cr}^{(2)}(k)=\frac{2^{k+1}( k-1)^{k+1}}{k^k}$. Then for the HC model (corresponding to the graph from Fig. 1) 
	\begin{itemize}
	\item for $\lambda\geq\lambda_{cr}^{(2)}$ there is exactly one 2-PNPGM $\mu_{0}$ which is translation invariant and corresponds to the solution $\left(A^*, A^* \right)$,
	\item for $0<\lambda<\lambda_{cr}^{(2)}$ there are exactly three 2-PNPGMs $\mu_{0}$, $\widehat{\mu}_{1}$ , $\widehat{\mu}_{2}$ corresponding to the solutions
 $\left(A^*, A^*\right), \ \left(A_1, A_2\right), \ \left(A_2, A_1\right),$
 \end{itemize}
where
$$A^*=\left(\dots,a^*,1,a^*,1,a^*,1,a^*,1,a^*,1,a^*,\dots\right), \ \ A_1=\left(\dots,a_1,1,a_1,1,a_1,1,a_1,1,a_1,1,a_1,\dots\right),$$
$$A_2=\left(\dots,a_2,1,a_2,1,a_2,1,a_2,1,a_2,1,a_2,\dots\right).$$}

\textbf{Case $q=4$.} By virtue of Corollary 6, the 4-periodic solutions (\ref{e18}) have the form:
\begin{center}
$z_j= \begin{cases}
1, \ \ \mbox{if} \ \ j\equiv0 \mod4, \\
a, \ \ \mbox{if} \ \ j\equiv1 \mod4, \\
\lambda_2,\ \ \mbox{if} \ \ j\equiv2 \mod4,\\
b, \ \ \mbox{if} \ \ j\equiv3 \mod4,
\end{cases}$ \ \ \ \
$\widetilde{z}_j= \begin{cases}
1, \ \ \mbox{if} \ \ j\equiv0 \mod4, \\
c, \ \ \mbox{if} \ \ j\equiv1 \mod4, \\
\lambda_2,\ \ \mbox{if} \ \ j\equiv2 \mod4,\\
d, \ \ \mbox{if} \ \ j\equiv3 \mod4.
\end{cases}$
\end{center}
In this case, we write the system of equations (\ref{e18}) as follows:
\begin{equation}\label{e21}
\begin{cases}
a=\lambda_1\Big({1+\lambda_2+c \over c+d}\Big)^k, \\
b=\lambda_3\Big({1+\lambda_2+d \over c+d}\Big)^k, \\
c=\lambda_1\Big({1+\lambda_2+a \over a+b}\Big)^k, \\
d=\lambda_3\Big({1+\lambda_2+b \over a+b}\Big)^k.
\end{cases}
\end{equation}
 In our further studies, it is very difficult to analyze the system of equations (\ref{e21}) in the general case.
 Therefore, in (\ref{e21}) we assume that $\lambda_{2i+1}=\lambda,$ $i=0,1$ and denote $\lambda_{2}=\gamma$. Then the system of equations (\ref{e21}) have form:

\begin{equation}\label{e22}
\begin{cases}
a=\lambda\Big({1+\gamma+c \over c+d}\Big)^k, \\
b=\lambda\Big({1+\gamma+d \over c+d}\Big)^k, \\
c=\lambda\Big({1+\gamma+a \over a+b}\Big)^k, \\
d=\lambda\Big({1+\gamma+b \over a+b}\Big)^k.
\end{cases}
\end{equation}

The next lemma is useful.

\textbf{Lemma 5}. \emph{For any} $k$ \emph{for the system of equations (\ref{e22}), the following statements are true}:

\emph{1.} $a=b$ \emph{if and only if} $c=d$.

\emph{2.} \emph{If} $a=c$, \emph{then} $b=d$.

\emph{3.} \emph{If} $b=d$, \emph{then} $a=c$.

Consider the mapping $W:R^{4}\rightarrow R^{4}$ defined as follows:
\begin{equation}\label{eq27}
\begin{cases}
a'=\lambda\Big({1+\gamma+c \over c+d}\Big)^k, \\
b'=\lambda\Big({1+\gamma+d \over c+d}\Big)^k, \\
c'=\lambda\Big({1+\gamma+a \over a+b}\Big)^k, \\
d'=\lambda\Big({1+\gamma+b \over a+b}\Big)^k.
\end{cases}
\end{equation}

Note that (\ref{eq27}) is the equation $x=W(x)$. To solve the system of equations (\ref{eq27}), we need to find the fixed points of the mapping $x'=W(x)$.

The following lemma is true.

\textbf{Lemma 6}. \emph{The following sets are invariant with respect to $W$ }:
$$I_{1}=\{x\in R^{4}: a=b=c=d\}, \ \ I_{2}=\{x\in R^{4}: a=c,~ b=d\},$$
$$I_{3}=\{x\in R^{4}: a=d,~ b=c\}, \ \ I_{4}=\{x\in R^{4}: a=b,~ c=d\}.$$

\textbf{Proof}. The invariance of $I_{1}$ and $I_{2}$ is obvious. Let us show the invariance of $I_{3}$ (the invariance of $I_{4}$ is proved similarly). It is clear that for any $x^{\ast}=(a^{\ast},b^{\ast},c^{\ast},d^{\ast})\in I_{3}$
$a^{ \ast}=d^{\ast},~ b^{\ast}=c^{\ast}$. We have
$$
\begin{cases}
a'=\lambda\Big({1+\gamma+c^{\ast} \over c^{\ast}+d^{\ast}}\Big)^k, \\
b'=\lambda\Big({1+\gamma+d^{\ast} \over c^{\ast}+d^{\ast}}\Big)^k, \\
c'=\lambda\Big({1+\gamma+a^{\ast} \over a^{\ast}+b^{\ast}}\Big)^k, \\
d'=\lambda\Big({1+\gamma+b^{\ast} \over a^{\ast}+b^{\ast}}\Big)^k,
\end{cases}
$$
i.e. $a'=d', \ b'=c'$, and hence $x'=W(x^{\ast})\in I_{3}$.
Lemma is proved.

\textbf{Remark 6}. Solutions of the system of equations (\ref{e22}) on the invariant sets $I_{1}$ and $I_2$ correspond to TINPGM. In particular, for $k\geq2$ the TINPGM corresponding to the collection of vectors on $I_{1}$ is unique, and on $I_2$ we obtain the assertion of Theorem 2.

\textbf{The case of $I_{3}$ and $k\geq2$.} The following assertion is true.

\textbf{Proposition 3}. Let $k\geq2$. Then for any values of $\lambda>0$ and $\gamma>0$, the system of equations (\ref{e22}) on the invariant set $I_{3}$ has a unique solution.

\textbf{Proof}. On $I_{3}$ the system of equations (\ref{e22}) has the form
\begin{equation}\label{e24}
\begin{cases}
a=\lambda\Big({1+\gamma+b \over a+b}\Big)^k, \\
b=\lambda\Big({1+\gamma+a \over a+b}\Big)^k.
\end{cases}
\end{equation}
In this system of equations, subtracting the second from the first equation, after some algebra we get

$$(a-b)\left[(a+b)^k+\lambda\left((1+\gamma+b)^{k-1}+\cdots+(1+\gamma+a)^{k-1 }\right)\right]=0,$$
i.e., $a=b$. Then in this case it follows from Lemma 2 that the system of equations (\ref{e24}), for any $\lambda>0$ and $\gamma>0$, has a unique solution of the form $(a^*, a^*)$. Proposition is proved.

\textbf{Case of $I_{4}$ and $k\geq2$}. In this case, from the system of equations (\ref{e24}) we get
\begin{equation}\label{e25}
a=g(c), \ \ \ c=g(a),
\end{equation}
where
$$g(x)=\lambda\left(\frac{1+\gamma+x}{2x}\right)^{k}.$$

The following proposition holds.

\textbf{Proposition 4}. Let $k\geq2$ and $\lambda_{cr}^{(3)}(k,\gamma)=\frac{2^{k}(\gamma+1)(k-1)^{k+1 }}{k^k}$. Then the system of equations (\ref{e25}) with $\lambda\geq\lambda_{cr}^{(3)}$ has exactly one solution of the form $\left(a^{\ast},a^{\ast} \right)$, and for $0<\lambda<\lambda_{cr}^{(3)}$ it has exactly three solutions of the form $\left(a^{\ast},a^{\ast}\right),$ $\left(a_1,a_2\right),$ $\left(a_2,a_1\right)$.

\textbf{Proof}. Let $h(x)=g(g(x))$. Let us show that the function $h(x)$ is $S$-shaped.
Indeed, the function $h(x)$ for $x>0$ is increasing and bounded, since
$$h'(x)=\frac{\lambda k^2(1+\gamma)^2\left(1+\gamma+g(x)\right)^{k-1}}{2^{k}x(x+\gamma+1)g^k(x)}>0, \ \ \ \   h(0)=\frac{\lambda}{2^k}>0,$$
$$\lim_{x\rightarrow\infty}{h(x)}=\lambda\left(\frac{2^{k}(1+\gamma)+\lambda}{2\lambda}\right)^k<\infty.$$
Moreover,
$$h''(x)=\frac{\lambda k^2(1+\gamma)^2\Big(g(x)\Big((k-1)(\gamma+1)-2x\Big)+\Big((k^2-1)(\gamma+1)-2x\Big)(1+\gamma)\Big)(1+\gamma+g(x)\Big)^{k-2}}
{2^{k}x^2(x+\gamma+1)^2g^k(x)}.$$

Let us show that the equation $h''(x)=0$ has exactly one solution. To do this, consider the function
$$r(x)=\lambda\left(\frac{1+\gamma+x}{2x}\right)^{k}\Big((k-1)(\gamma+1)-2x\Big)+\Big((k^2-1)(\gamma+1)-2x\Big)(1+\gamma).$$

Calculate the derivative $r'(x)$:
$$r'(x)=-\frac{g(x)}{x^2+(1+\gamma)x}\Big(2x^2-2(k-1)(1+\gamma)x+k(1+\gamma)^2(k-1)\Big)-2(\gamma+1),$$
where $g(x)=\lambda\left(\frac{1+\gamma+x}{2x}\right)^{k}.$
It is easy to see that
$$2x^2-2(k-1)(1+\gamma)x+k(1+\gamma)^2(k-1)>0,$$
where $k\geq2$.

Hence $r'(x)<0$, i.e. function $r(x)$ is decreasing.
On the other hand, for $k\geq2$
$$r\left(\frac12\right)=\lambda\Big(\frac{3+2\gamma}{2}\Big)^k\Big((k-1)\gamma+k-2\Big)+\Big((k^2-1)\gamma+k^2-2\Big)(1+\gamma)>0, $$
$$r\left(\frac{(k^2-1)(\gamma+1)}2\right)=-\lambda k(k-1)(\gamma+1)\Big(\frac{2+2\gamma+(k^2-1)(\gamma+1)}{2(k^2-1)(\gamma+1)}\Big)^k<0.$$
Therefore, the equation $r(x)=0$ and hence the equation $h''(x)=0$ has a unique root for $x>0$. Hence the function $h(x)$ is $S$-shaped.

It is easy to see that since $h'(x_0)=(g'(x_0))^2$, then $x_0$ is the only fixed point of the function $h(x)$ if and only if $g'(x_0)\geq-1$.

Let $x_0$ be the unique solution of the equation $g(x)=x$. Calculate the derivative $g'(x_0)$:
$$g'(x)=-\frac{\lambda k(1+\gamma+x)^{k-1}(\gamma+1)}{2^{k}x^{k+1}} \ \ \Rightarrow \ \ g'(x_0)=-\frac{k(1+\gamma)}{1+\gamma+x_0}.$$
After solving the inequality $g'(x_0)<-1$, we have $x_0<(k-1)(\gamma+1)$. Then from $x_0=\lambda\left(\frac{1+\gamma+x_0}{2x_0}\right)^{k}$ we get that the equation $h(x)=x$ has only one fixed point at $\lambda\geq\frac{2^{k}(\gamma+1)(k-1)^{k+1}}{k^k}=\lambda_{cr}^{(3)}(k ,\gamma).$ From the inequality $g'(x_0)<-1$ we obtain $\lambda<\lambda_{cr}^{(3)}$. Then by Kesten's lemma it follows that the function $h(x)$ has at least three fixed points. On the other hand, under this condition, by the property of an $S$-shaped function, we have at most three fixed points for $h(x)$. Hence, for $0<\lambda<\lambda_{cr}^{(3)}$ there are exactly three fixed points of the equation $h(x)=x$.
Therefore, the system of equations (\ref{e25})
for $\lambda\geq\lambda_{cr}^{(3)}$ has a unique solution of the form $(a^{\ast},a^{\ast})$, and for $\lambda<\lambda_{cr }^{(3)}$ it has two positive solutions $\Big(a_1,a_2\Big)$ and $\Big( a_2,a_1\Big)$. The proposition is proved.

Thus, for $k\geq2$ and $0<\lambda<\lambda_{cr}^{(3)}$, the coordinates of the vector $(z, \widetilde{z})$, which is the solution of the system of equations (\ref{e18}), have the following forms:
$$
\begin{cases}
z_{4i}=\widetilde{z}_{4i}=1, \ \ \ i\in \mathbb{Z}, \\
z_{4i+2}=\widetilde{z}_{4i+2}=\gamma, \ \ \ i\in \mathbb{Z}, \\
z_{2i+1}=a_1, \ \ \ i\in \mathbb{Z},\\
\widetilde{z}_{2i+1}=a_2, \ \ \ i\in \mathbb{Z}
\end{cases}
$$
and
$$
\begin{cases}
z_{4i}=\widetilde{z}_{4i}=1, \ \ \ i\in \mathbb{Z}, \\
z_{4i+2}=\widetilde{z}_{4i+2}=\gamma, \ \ \ i\in \mathbb{Z}, \\
z_{2i+1}=a_2, \ \ \ i\in \mathbb{Z},\\
\widetilde{z}_{2i+1}=a_1, \ \ \ i\in \mathbb{Z}.
\end{cases}
$$

Hence, it is clear that $\sum_{i\in \mathbb Z} z_{i}=+\infty$
and $\sum_{i\in \mathbb Z} \widetilde{z}_{i}=+\infty,$ i.e. the series obtained from the sequence of solutions diverge.

By Propositions 3 and 4, the following theorem holds.

\textbf{Theorem 4.} \textit{Let $k\geq2$ and $q=4$. Then the following statement is true for the HC model (corresponding to the graph from Fig. 1):}

1. \textit{For $\lambda>0$ and $\gamma>0$, there is exactly one 4-PNPGM on $I_3$. Moreover, this measure coincides with the only TINPGM.}

2. \textit{Let $\lambda_{cr}^{(3)}(k,\gamma)=\frac{2^{k}(\gamma+1)(k-1)^{k+1} }{k^k}$. Then on $I_4$ for
$\lambda\geq\lambda_{cr}^{(3)}$ and $\gamma>0$, there is exactly one 4-PNPGM $\mu_{0}$ which is translation invariant and corresponds to the solution $\left (A^*, A^*\right)$, and for $0<\lambda<\lambda_{cr}^{(3)}$ and $\gamma>0$, there are exactly three 4-PNPGM $\mu_{ 0}$, $\widehat{\mu}_{1}$, $\widehat{\mu}_{2}$ corresponding to the solutions $\left(A^*, A^*\right), \ \left (A_1, A_2\right), \ \left(A_2, A_1\right),$
where
$$A^*=\left(\dots ,a^*,1,a^*,\gamma,a^*,1,a^*,\gamma,a^*,1,a^*,\dots\right),\ \ A_1=\left(\dots,a_1,1,a_1,\gamma,a_1,1,a_1,\gamma,a_1,1,a_1,\dots\right),$$
$$A_2=\left(\dots,a_2,1,a_2,\gamma,a_2,1,a_2,\gamma,a_2,1,a_2,\dots\right).$$}

\textbf{Acknowledgements.} This work was supported by the fundamental project
(number: F-FA-2021-425) of the Ministry of Innovative Development of the Republic of Uzbekistan.
Rozikov thanks Institut des Hautes \'Etudes Scientifiques (IHES), Bures-sur-Yvette, France and the IMU-CDC for support of his visit to IHES.\\

\section*{Statements and Declarations}

{\bf	Conflict of interest statement:} 
On behalf of all authors, the corresponding author (U.A.Rozikov) states that there is no conflict of interest.

\section*{Competing Interests and Funding}
This work was supported by the fundamental project
(number: F-FA-2021-425) of the Ministry of Innovative Development of the Republic of Uzbekistan.

\section*{Data availability statements}
The datasets generated during and/or analysed during the current study are available from the corresponding author on reasonable request.


\begin{thebibliography}{99}
	
\bibitem{RM} R.M. Khakimov, M.T. Makhammadaliev, Uniqueness and nonuniqueness conditions for weakly periodic Gibbs measures for the Hard-Core model, {\it Theor. Math. Phys.}, 2020, V.204, No.2, p. 1059-1078.

\bibitem{R} U.A. Rozikov, \emph{Gibbs measures on Cayley trees}, World Scientific. 2013.

\bibitem{BR} L.V. Bogachev, U.A. Rozikov, On the uniqueness of Gibbs measure in the Potts model on a Cayley tree with external field, \emph{J. Stat. Mech. Theory Exp.}, 2019, no. 7, 073205, 76 pp.

\bibitem{GR}  N.N. Ganikhodjaev, U. A. Rozikov, The Potts Model with Countable Set of Spin Values on a Cayley Tree, \emph{Letters in Mathematical Physics}, 75, (2006), 99-109.

\bibitem{G}  N.N. Ganikhodjaev, Limiting Gibbs measures of Potts model
with countable set of spin values, \emph{J. Math. Anal. Appl.}, 336 (2007), 693-703.

\bibitem{Z}   Ye Zichun, Models of gradient type with sub-quadratic actions,  \emph{J. Math. Phys.}, 60, 073304 (2019). 

\bibitem{HKR}  F. Henning, C. K\"{u}lske, A. Le Ny, U.A. Rozikov, Gradient gibbs measures for the SOS-model with countable values on a Cayley tree, \emph{Electron. J. Probab.}, vol.24, 2019.  DOI: 10.1214/19-EJP364.

\bibitem{HK}   F. Henning, C. K\"{u}lske, Coexistence of localized Gibbs measures and delocalized gradient Gibbs measures on trees, \emph{Ann. Appl. Probab.}, 31 (5), 2284-2310, (2021). 

\bibitem{B}   S. Buchholz, Phase transitions for a class of gradient fields, \emph{Probability Theory and Related Fields}, (2021) 179, 969-1022.

\bibitem{HK1}  F. Henning, C. K\"{u}lske, Existence of gradient Gibbs measures on regular
trees which are not translation invariant,  arXiv:2102.11899v2.

\bibitem{RM1} R.M. Khakimov, M.T. Makhammadaliev, Nonprobability Gibbs measures for the HC
model with a countable set of spin values for a Wand type graph on a Cayley tree, {\it Theor. Math. Phys.} 2022, V.212, No.3, p. 1259-1275.
	
\bibitem{AI} E. Aghion, D.A. Kessler, E. Barkai: From non-normalizable Boltzmann-Gibbs statistics to infinite-ergodic theory, {\em Phys. Rev. Lett}. \textbf{122}   (2019), 010601 (5 pages).

\bibitem{bw} G. Brightwell, P. Winkler, Graph homomorphisms and phase transitions, \emph{J.Combin. Theory Ser.B.}, 1999. {\bf 77}, P. 221-262.

\bibitem{Kel}	D. Galvin, F. Martinelli, K. Ramanan, P. Tetali, The multi-state Hard Core model on a regular tree,  \emph{SIAM Journal on Discrete Mathematics},  2011. \textbf{25}(2), 894-915.

\bibitem{K} H. Kesten, Quadratic Transformations: A Model for Population Growth. I, \textit{Adv. Appl. Probab.}, \textbf{2}:1 (1970), 1-82.


\end{thebibliography}
\end{document}